\title{Stochastic Stability for  Flows
with Smooth Invariant Measures}

\documentclass[11pt]{article}
\usepackage{amsmath}
\usepackage{epsfig}
\usepackage{graphics}
\usepackage{amsfonts}
\usepackage{amssymb}

\small\normalsize

\newtheorem{Thm}{Theorem}[section]  
\newtheorem{Prop}[Thm]{Proposition}

\numberwithin{equation}{section}

\newcommand{\D}[1]{{\mathbb#1}}

\newcommand{\RR}{{\D{R}}}

\newcommand{\TT}{{\D{T}}}

\newtheorem{defn}{Definition}[section]

\newtheorem{thm}[defn]{Theorem}

\newtheorem{prop}[defn]{Proposition}

\author{
Sergiu Aizicovici and Todd Young \\
Department of Mathematics\\
Ohio University, Ohio, USA\\
{\small azicovs@ohio.edu}\\
{\small youngt@ohio.edu}
}

\begin{document}

\maketitle

\begin{abstract}
We study the notion of stochastic stability with respect to diffusive perturbations
for flows with smooth invariant measures. We investigate the question fully for non-singular flows on
the circle. We also show that  volume-preserving flows are
stochastically stable with respect to perturbations that are
associated with homogeneous diffusions.
\end{abstract}

Published in {\it Libertas Matematica} {\bf 30} (2010), 71-79.

\section{Introduction}

The notion of stochastic stability
of the invariant measures of a continuous dynamical system was formulated as early as
1933 by Pontryagin, Andronov and Vitt \cite{PAV} and Bernstein \cite{Ber}.
They solved the problem completely for flows on the real line that possess a globally attracting
set of equilibrium points. Such flows of course are the gradients of some potential function.
Freidlin and Wentzell \cite{FW} extended the classical work to $\RR^n$, again with globally
attracting equilibrium points, by defining dynamically a quasi-potential associated with
the flow. In the case of a unique equilibrium, stochastic stability is proved and
a formula for the perturbed invariant measure was produced in terms of the quasi-potential.
Stochastic stability for flow on a compact manifold was discussed in
\cite{Ki}, \cite{Ki91} and extended to the case that the maximal invariant set of the flow consists of
a finite number of hyperbolic basic sets.
Since then, it seems that very little attention has been given to this
problem.

In the meantime, stochastic stability of invariant measures has been a topic of
continuing interest and surprising difficulty for discrete time dynamical systems theory.
The main directions where progress has been made are: structurally stable systems \cite{Ki90}, hyperbolic systems \cite{Ki86}, \cite{Ki90}, \cite{Yo86} and unimodal maps \cite{BY92}, \cite{BV96}.

In this note we study the question of stochastic stability for the
(overlooked) case of non-singular flow on the circle and for a
volume preserving flow on a manifold.

We will let $h: M\to TM$ be a smooth vector field on a compact Riemannian manifold $M$
with the associated differential equation in local coordinates:
\begin{equation}\label{ode}
 \frac{dx}{dt} = h(x).
\end{equation}
By standard results, $h$ generates a smooth, global solution flow
$\Phi^t:\RR \times M \rightarrow M$.
All equations in this note will be assumed to
be written in local coordinates.

We will denote by $d$ the distance on $M$, by $\mathcal{B}$ the Borel
$\sigma$-algebra on $M$, by $m$ the normalized Riemannian volume (Lebesgue measure)
on $M$, and by $\mathcal{M}$ the space of Borel probability measures on $M$ with
the weak topology. We will use $\rightharpoonup$ to denote weak convergence in
$\mathcal{M}$.

Consider a vector field $h$ for which there is a unique physical,
ergodic measure $\mu_0 \in \mathcal{M}$.
An ergodic measure $\mu_0$ is {\em physical} if for $m$-a.e.~$x \in M$,
$\nu_{x,T} := \frac{1}{T} \int_{0}^{T} \delta_{\Phi^t(x)} \rightharpoonup \mu_0$
as $n \rightarrow \infty$,
i.e.~for all $\phi\in C^0(M)$ we have
$\frac{1}{T} \int_{0}^{t} \phi(\Phi^t(x)) \rightarrow \int \phi \, d\mu_0,$.
We will also assume that $\mu_0$ is absolutely continuous
with respect to the Lebesgue measure on the manifold and has a smooth density function
$\rho_0$.

By {\em stochastic stability} we mean stability of $\mu_0$ under small stochastic perturbations
of $\phi^t$, which we now define. Consider the stochastic differential equation
\begin{equation}\label{sde}
 dx = h(x) \, dt + \sqrt{\epsilon} \gamma(x) \circ dW,
\end{equation}
which we interpret as a Stratanovich integral equation on the tangent space of the manifold $M$.
Equation (\ref{sde}) is associated with a diffusion process on $M$.
For general background on diffusions on a manifold see for instance
\cite[Chp. V]{IW}, \cite{El} or \cite{Fr04}.
We call (\ref{sde}) a small stochastic perturbation of (\ref{ode}).
A class of perturbations refers to a collection of equations
of the form (\ref{sde}) with $\gamma(x)$ chosen from some class
of nonsingular matrix fields.
By stochastic stability with respect to a class  of
small stochastic perturbations, we mean that each flow defined by (\ref{sde}) within
the class has a unique ergodic stationary measure $\mu_\epsilon$
and $\mu_\epsilon \rightharpoonup \mu_0$ as $\epsilon \rightarrow 0$.

Note that any weak limit point $\mu$ of the set $\{ \mu_\epsilon \}$ as
$\epsilon \rightarrow 0$ is called a {\em zero-noise limit} measure for $h$.
It is well-known that all zero-noise limit measures are invariant under $h$ \cite{Ki}.
See also \cite{CY05} for a recent related work using zero-noise limits.

We note here that no counterexamples to stochastically stability for
physical ergodic measures of dynamical systems, except one in which
the ``random'' perturbations have a distinctly non-random character \cite{AT}.

By {\em standard diffusion} we  mean for each $x$, $\gamma(x)$ is a multiple of the identity matrix.

Also note that we could consider larger classes of perturbations that include both
stochastic and deterministic perturbations by adding a term
$\epsilon g \, dt$ to the right hand side of (\ref{sde}), where $g$ is
selected from some class of vector fields.

We can associate with (\ref{ode}) and (\ref{sde}) a generator
\begin{equation}\label{generator}
\mathcal{L}_\epsilon = h \cdot \nabla +  \frac{\epsilon}{2} \Gamma : \nabla \nabla
\end{equation}
where $\Gamma(z) = \gamma(x) \gamma(x)^T$ and ``$:$'' denotes the double inner product
of matrices (the inner product consistent with the Frobenius norm).
The formal $L^2$ adjoint of $\mathcal{L}$ is
\begin{equation}\label{adjoint}
 \mathcal{L}^*_\epsilon = - \nabla \cdot h
      +  \frac{\epsilon}{2} \nabla \cdot \nabla \cdot \Gamma.
\end{equation}
Note that $\Gamma(x)$ is a symmetric positive definite matrix field.
For standard diffusion, $\Gamma(x)$ is the identity matrix and $\mathcal{L}_\epsilon$
and $\mathcal{L}^*_\epsilon$ reduce to:
\begin{equation}\label{generator_standard}
\mathcal{L}_\epsilon = h \cdot \nabla + \frac{\epsilon}{2}  \Delta
\qquad \textrm{and} \qquad
\mathcal{L}^*_\epsilon = - \nabla \cdot h
      + \frac{\epsilon}{2} \Delta.
\end{equation}
Here $\Delta$ is the Laplace-Beltrami operator on $M$.

Probability density functions evolve under (\ref{sde}) via the
Fokker-Planck equation
\begin{equation}\label{Fokker-Plank}
\frac{\partial \rho}{\partial t} = \mathcal{L}^*_\epsilon \rho.
\end{equation}
(For $\epsilon =0$ this is usually called the Liouville equation.)
Suppose that $\mu_0$ is an invariant measure for (\ref{ode})  and the
density $\rho_0$ of $\mu_0$ is a $C^2$ function.
Then $\mathcal{L}^*_0 \rho_0 =0$.
Further, if $\rho_\epsilon$ is the density of a stationary measure $\mu_\epsilon$ for
(\ref{sde}), then
\begin{equation}\label{stationary}
   \mathcal{L}^*_\epsilon \rho_\epsilon =0.
\end{equation}
Since $M$ is compact and $\Gamma$ is positive definite,
the density $\rho_\epsilon$ for the perturbed system (\ref{sde})
exists \cite{Do} and it is unique, nonzero and smooth \cite[Prop.~5.4.5]{IW}.

Now we wish to treat the behaviour of $\rho_\epsilon$ as a perturbation problem.
Define $r_\epsilon(x) = \rho_\epsilon(x) - \rho_0(x)$ so that
\begin{equation}\label{ansatz}
 \rho_\epsilon = \rho_0 + r_\epsilon.
\end{equation}
It follows that  $r_\epsilon$ exists, is unique and is smooth because
$\rho_\epsilon$ exists and is unique and smooth, and, $\rho_0$ is unique and smooth
by assumption. Since $\int \rho_\epsilon = 1$, then any solution
$r_\epsilon$ must satisfy: $\int r_\epsilon = 0$.

Substitution into (\ref{stationary}) gives a perturbed elliptic equation with a constraint:
\begin{equation}\label{r_pde}
\frac{\epsilon}{2} \, \nabla \cdot \nabla \cdot (\Gamma r_\epsilon)
    - \nabla \cdot (hr_\epsilon)
      =  - \frac{\epsilon}{2} \,  \nabla \cdot \nabla \cdot (\Gamma \rho_0),
          \quad \textrm{ subject to } \int r_\epsilon = 0.
\end{equation}
In the case of standard diffusion, this simplifies to:
\begin{equation}\label{r_pde_standard}
\frac{\epsilon}{2} \, \Delta r_\epsilon =  \nabla \cdot (hr_\epsilon)
        - \frac{\epsilon}{2} \, \Delta \rho_0,
          \quad \textrm{ subject to } \int r_\epsilon = 0.
\end{equation}
To complete this approach we need only show that $r_\epsilon$ goes to
zero as $\epsilon$ goes to zero.

{\bf Note:}\\
Instead of (\ref{ansatz}) we could suppose that $\rho_\epsilon$ has the form
\begin{equation}
    \rho_\epsilon(x) = \rho_0(x) + \epsilon u(x) + r_\epsilon(x)
\end{equation}
for some smooth function $u(x)$ and remainder $r_\epsilon(x)$ of order $o(\epsilon)$.
However, this is a stronger assumption than we need and may impose unnecessary
constraints. In one dimension one can solve directly for $u$, but in higher dimensions
a solution seems unlikely.


\section{Two special cases}

In this section we point out that stochastic stability with respect to
standard diffusion follows easily from known results in two extreme
cases, volume preserving flows and gradient flows.

First suppose that $h$ is volume preserving and the perturbation
is by standard diffusion. According to \cite[Theorem~5.4.6]{IW}
the measure $c m$ ($c$ constant) is an invariant measure of a diffusion
generated by $\Delta + h$ if and only if $h$ is volume preserving.
Rescaling (\ref{sde}) by $1/\epsilon$ implies that the volume itself
is the invariant measure of the perturbed process for any $\epsilon >0$.
This gives us the following:
\begin{Prop}
 If $h$ is volume-preserving and the volume, $m$, is ergodic, then $m$
is stochastically stable under standard diffusion.
\end{Prop}
We will extend this result slightly to the case of homogeneous diffusions
in \S~4.

Secondly, consider the case that $h$ is a gradient flow, i.e., there exists
a smooth, real-valued function $H$ on $M$ such that $h = - \textrm{grad} \, H$.
This case is the direct generalization of the early results \cite{PAV}, \cite{Ber} and
a special case of result in \cite{FW}.
In \cite[Theorem~5.4.6]{IW} we find that if $h = -\textrm{grad} \, H$ then the diffusion
generated by $\Delta + h$ has as its invariant measures multiples of
$\exp(-2H(x)) \, dx$. Rescaling by $1/\epsilon$ implies that the
invariant measure for \ref{sde} is:
\begin{equation}\label{gradmeas}
     c_\epsilon  e^{-2H(x) / \epsilon} \, dx.
\end{equation}
This measure converges as $\epsilon \rightarrow 0$ to a measure
supported on the set of minimum points of the function $H(x)$.
In the simplest case we have:
\begin{Prop}
 If $h$ is a gradient flow with potential $H$, and $H$ has
a unique minimum point, then the delta measure
on this point is stochastically stable.
\end{Prop}
If there is not a unique minimum, then the delta measure
at each of the minimum points is ergodic and any linear
combination of them is (nonergodic) invariant. Freidlin
has studied the asymptotics of transitions between the
local minima as $\epsilon \rightarrow 0$ \cite{Fre}.


\section{Flow on a Circle}

Suppose $M = \TT^1$ and $h$ is smooth and nonzero. Then $h$ generates
a smooth flow and possess a smooth ergodic measure supported on the entire
circle, given by the density $\rho_0(x) = c/h(x)$.

On the circle, (\ref{r_pde}) becomes:
\begin{equation}\label{r_pde_1d}
\frac{\epsilon}{2} \,  (\Gamma r_\epsilon)'' -  (h r_\epsilon)'
     =   - \frac{\epsilon}{2} \, (\Gamma \rho_0)''
          \quad \textrm{ subject to } \int r_\epsilon = 0.
\end{equation}

For $\epsilon =0$, note that there is a one dimensional space of solutions of
(\ref{r_pde_1d}); $r_0 = c/h(x)$ is a solution for any constant $c$.
If we restrict to $\int r_0 = 0$, then $r_0 \equiv 0$ is the unique solution.
As stated before, the existence, uniqueness and smoothness of the solution $r_\epsilon$
of (\ref{r_pde_1d}) follow from the assumptions on $\rho_0$ and the existence, uniqueness and
smoothness of $\rho_\epsilon$. However, there is also an elementary proof of these facts in this case.

Assume that $h$, $\rho_0$ and $\Gamma$ are
all at least $C^2$ smooth on the circle.
Then $r_\epsilon$ is a solution of (\ref{r_pde_1d})
if and only if it is a solution of
\begin{equation}\label{r_pde_1d_reduced}
 \frac{\epsilon}{2} \,  (\Gamma r_\epsilon)' -  h r_\epsilon
     =   - \frac{\epsilon}{2} \, (\Gamma \rho_0)' + C_\epsilon,
\end{equation}
for some constant $C_\epsilon$.

\begin{prop}
If $r_\epsilon$ is a solution  of (\ref{r_pde_1d}),
then $r_\epsilon$ converges to zero in the $L^2$ norm as $\epsilon \rightarrow 0$.
\end{prop}
\noindent
{\sf Proof:}\\
We will assume without loss of generality that $h(x)>0$ in local coordinates.
Let $r_\epsilon$ be a solution of (\ref{r_pde_1d}) subject to $\int r =0$.
Multiplying equation (\ref{r_pde_1d_reduced}) by $r_\epsilon$ and integrating
over $\TT^1$ we have:
$$
   \frac{\epsilon}{2} \,  \int (\Gamma r_\epsilon)' r_\epsilon \, dx
          - \int h r_\epsilon^2  \, dx
     =   - \frac{\epsilon}{2} \, \int (\Gamma \rho_0)' r_\epsilon \, dx
         + \int C r_\epsilon \, dx.
$$
Here and subsequently, $\int$ will mean integration over the circle, $\TT^1$.
Note first that the final integral  $\int C_\epsilon r_\epsilon \, dx$ is zero by constraint.
Rearranging and successively integrating by parts we obtain:
\begin{equation*}
\begin{split}
 \int h r_\epsilon^2  \, dx
     &= - \frac{\epsilon}{2} \,  \int \Gamma r_\epsilon r_\epsilon' \, dx
         + \frac{\epsilon}{2} \, \int (\Gamma \rho_0)' r_\epsilon \, dx\\
     &= - \frac{\epsilon}{4} \,  \int \Gamma \frac{d}{dx} r_\epsilon^2  \, dx
          + \frac{\epsilon}{2} \, \int (\Gamma \rho_0)' r_\epsilon \, dx \\
     &=  \frac{\epsilon}{4} \,  \int \Gamma' r_\epsilon^2  \, dx
         + \frac{\epsilon}{2} \, \int (\Gamma \rho_0)' r_\epsilon \, dx.
\end{split}
\end{equation*}
Thus we have
$$
   \int (h - \frac{\epsilon}{4} \Gamma') r_\epsilon^2  \, dx
        =  \frac{\epsilon}{2} \, \int (\Gamma \rho_0)' r_\epsilon \, dx.
$$
Denote $\alpha = \min h$. In the case that $\Gamma$ is not identically constant,
note that we have $\max(\Gamma') > 0$ and so
\begin{equation}\label{hGamma}
h - \frac{\epsilon}{4} \Gamma' > \alpha/2
\end{equation}
provided that $\epsilon < 2 \alpha / \max(\Gamma')$. (If $\Gamma$ is identically constant,
then (\ref{hGamma}) holds for any $\epsilon$.)
Therefore,
$$
 \frac{\alpha}{2} \int  r_\epsilon^2  \, dx
        \le \frac{\epsilon}{2} \, \int |(\Gamma \rho_0)' r_\epsilon| \, dx.
$$
Applying the Schwartz inequality to the last integral we have
$$
     \frac{\alpha}{2}  \| r_\epsilon\|_2^2
        \le \frac{\epsilon}{2} \|(\Gamma \rho_0)'\|_2 \|r_\epsilon\|_2 .
$$
or
\begin{equation}\label{l2r}
    \| r_\epsilon\|_2
        \le \frac{\epsilon}{\alpha} \|(\Gamma \rho_0)'\|_2 .
\end{equation}
Thus we have in fact that $\|r_\epsilon\|_2 = O(\epsilon)$.
\hfill $\Box$

Next we show that $r'_\epsilon = O(\epsilon)$ in the $L^2$ norm.
Multiplying equation (\ref{r_pde_1d}) by $r'_\epsilon$ and integrating
over $\TT^1$ we have:
$$
   \frac{\epsilon}{2} \,  \int (\Gamma r_\epsilon)'' r_\epsilon' \, dx
          - \int (h r_\epsilon)' r_\epsilon' \, dx
     =   - \frac{\epsilon}{2} \, \int (\Gamma \rho_0)'' r_\epsilon' \, dx.
$$
Expanding the derivatives and integrating by parts we obtain:
$$
    \int (h - \frac{3 \epsilon}{4} \Gamma') (r_\epsilon')^2  \, dx
        = - \int (h' - \frac{\epsilon}{2} \Gamma'') r_\epsilon r_\epsilon'  \, dx
            +  \frac{\epsilon}{2} \, \int (\Gamma \rho_0)'' r_\epsilon' \, dx.
$$
Note that for $\epsilon < 2 \alpha / ( 3 \max \Gamma')$ (or any $\epsilon$ in the case $\Gamma' \equiv 0$)
we have $h - \frac{3 \epsilon}{4} \Gamma' > \alpha/2$, and so
$$
 \frac{\alpha}{2} \int  (r_\epsilon')^2  \, dx
        \le \int |h' - \frac{\epsilon}{2} \Gamma''| |r_\epsilon r_\epsilon'|  \, dx
            +  \frac{\epsilon}{2} \, \int |(\Gamma \rho_0)'' r_\epsilon'| \, dx.
$$
Let
$$
\beta = \max( |h'| + \frac{\alpha}{3 \max \Gamma'}| \Gamma''|).
$$
(In the case $\Gamma' \equiv 0$, we may take $\beta = \max( h')$.
Then using the Schwartz inequality, we obtain:
\begin{equation}\label{l2rprime}
    \| r_\epsilon' \|_2^2
        \le \frac{2\beta}{\alpha} \|r_\epsilon\|_2 \|r_\epsilon'\|_2
            + \frac{\epsilon}{\alpha} \|(\Gamma \rho_0)''\|_2 \| r'_\epsilon \|_2.
\end{equation}
Using  (\ref{l2r}) we have:
$$
\| r_\epsilon' \|_2
        \le \epsilon \left(
      \frac{2\beta}{\alpha^2}   \|(\Gamma \rho_0)'\|_2
                  + \frac{1}{\alpha}  \|(\Gamma \rho_0)''\|_2 \right).
$$
Thus  we have:
\begin{prop}
As $\epsilon$ goes to zero, $\| r_\epsilon' \|_2$ is of order epsilon.
\end{prop}

Now, by the Poincar\'{e}-Wirtinger
inequality in one dimension \cite[p.~146.]{Bre}, the $L^\infty$ norm of $r_\epsilon$ is also of order $O(\epsilon)$. (An elementary proof of this also exists in this case.)
Since $r_\epsilon$ is smooth and the circle is compact, the following proposition holds.
\begin{prop}
The solution $r_\epsilon$ of (\ref{r_pde_1d}) converges to zero uniformly
as $\epsilon \rightarrow 0$.
\end{prop}

In terms of Stochastic Stability we have shown the following result.
\begin{thm}
Let $h \in C^1$ be a nonsingular flow on the circle with an
absolutely continuous invariant measure $\mu_0$
with density $\rho_0 \in C^2$. Then $\mu_0$ is stochastically stable with respect to
any perturbations in the class $C^2$.
\end{thm}


\section{Volume-Preserving Flows}

Next consider smooth vector fields $h$ that preserve volume, thus $\nabla \cdot h = 0$.
If volume is also  ergodic for the flow defined by $h$, then the constant multiples of the volume
are the only absolutely continuous invariant measures.

In this case (\ref{r_pde}) becomes:
\begin{equation}\label{r_pde_vp}
\frac{\epsilon}{2} \, \nabla \cdot \nabla \cdot (\Gamma r_\epsilon)
    - \nabla \cdot (hr_\epsilon)
      =  - \frac{\epsilon c}{2} \,  \nabla \cdot \nabla \cdot \Gamma,
          \quad \textrm{ subject to } \int r_\epsilon = 0.
\end{equation}

We note that $\nabla \cdot h = 0$ implies:
$$
    \nabla \cdot (hr_\epsilon) = h \cdot \nabla r_\epsilon
$$

If $\Gamma$ is a multiple of the identity, or, if $\Gamma$ is constant with respect to $x$ (homogeneous),
then (\ref{r_pde_vp}) simplifies to:
\begin{equation}\label{r_pde_vp_stand}
    \nabla \cdot (hr_\epsilon) = \frac{\epsilon}{2} \, \Delta r_\epsilon,
          \quad \textrm{ subject to } \int r_\epsilon = 0.
\end{equation}
\begin{prop}
 The only smooth solution of (\ref{r_pde_vp_stand}) is $r_\epsilon \equiv 0$.
\end{prop}

\noindent
{\sf Proof:}
Multiplying both sides of the equation by $r_\epsilon$ and integrating we have:
\begin{equation}\label{hr1}
\begin{split}
  \int (h \cdot \nabla r_\epsilon) r_\epsilon
      &= \frac{\epsilon}{2} \int r_\epsilon \Delta r_\epsilon \\
   \int (h r_\epsilon) \cdot \nabla r_\epsilon  &= -  \frac{\epsilon}{2} \int | \nabla r_\epsilon |^2
\end{split}
\end{equation}

On the other hand, $r_\epsilon$ is a smooth solution of (\ref{r_pde_vp_stand}) if and only
if there is a smooth divergence-free vector field $v_\epsilon$ such that
$$
    h r_\epsilon = \frac{\epsilon}{2} \, \nabla r_\epsilon + v_\epsilon.
$$
We may take the inner product of both sides of the equation with $\nabla r_\epsilon$ and integrate
to obtain:
\begin{equation}\label{hr2}
\begin{split}
   \int (h r_\epsilon) \cdot \nabla r_\epsilon
            &= \frac{\epsilon}{2} \int | \nabla r_\epsilon |^2
                +  \frac{\epsilon}{2} \int v_\epsilon \cdot \nabla r_\epsilon \\
            &= \frac{\epsilon}{2} \int | \nabla r_\epsilon |^2
                - \frac{\epsilon}{2} \int \nabla v_\epsilon \cdot r_\epsilon \\
            &= \frac{\epsilon}{2} \int | \nabla r_\epsilon |^2
\end{split}
\end{equation}
We can conclude from (\ref{hr1}) and (\ref{hr2}) that both
$\int (h r_\epsilon) \cdot \nabla r_\epsilon$ and $\int | \nabla r_\epsilon |^2$ are zero.
This along with the constraint $\int r_\epsilon = 0$ implies the result. \hfill $\Box$

We note here that there are other simple proofs that $\int (h r_\epsilon) \cdot \nabla r_\epsilon =0$,
including a dynamical argument. To see this, note that
$(h r_\epsilon) \cdot \nabla r_\epsilon = (h \cdot \nabla r_\epsilon) r_\epsilon$ is in fact the
directional derivative of $r^2$ along the vector field $h$. If we integrate this expression
along a trajectory, ergodicity and the assumption that $r_\epsilon$ is bounded imply the claim.

The implications for stochastic stability are the following:
\begin{Thm}
If $h$ preserves volume, then the volume
is the unique zero-noise limit under any homogeneous diffusion.
If volume is ergodic, then it is stochastically stable with
respect to homogeneous diffusion.
\end{Thm}


\section{Discussion of non-smooth invariant densities}

In general, it is rare for $\mu_0$ to have a smooth
density, rather, it is usually supported on an attractor which
has dimension less than that of the ambient space. Further complications
include that the physical measure $\mu_0$ will not be the unique ergodic measure supported on
the attractor.

If the method could be generalized it might be as follows. The Liouville equation
must be considered as acting on distributions. The densities for the SDE
(\ref{sde}) however
remain smooth. With these considerations, (\ref{ansatz})
cannot be the right assumption. Rather, we need to show something of the form
$$
   \rho_\epsilon = A_\epsilon \rho_0,
$$
where $A_\epsilon$ is a smoothing operator that acts on distributions and
approaches the identity as $\epsilon \rightarrow 0$. It is natural
to assume then is that $A_\epsilon$ is itself generated by a diffusion, i.e.
$$
  A_\epsilon = e^{\epsilon L^*}
$$
where $L = \hat{\Gamma}(x) \Delta$ and $\hat{\Gamma}$ must be determined.
In this circumstance, with standard diffusion as the class of perturbation,
$\mathcal{L}^*_\epsilon \rho_\epsilon = 0$ becomes:
$$
   \nabla \cdot h (A_\epsilon \rho_0) = \frac{\epsilon^2}{2} \Delta (A_\epsilon \rho_0),
$$
and the problem is to find an appropriate $\hat{\Gamma}(x)$.

\smallskip
{\bf Acknowledgment:} The authors wish to thank Martin Hairer for
pointing out important references and making several helpful comments.

\end{document}